%Paper: 9202005
%From: Guenter Ziegler <ziegler@ml.kva.se>
%Date: Tue, 11 Feb 92 08:35:04 +0100

%
% prepared with plain-TeX, C-Version 3.14
%
\magnification=\magstep1
\baselineskip=13pt
\font\typc=cmbx10 scaled \magstep1  % for unnumbered section titles
\font\type=cmbx10 scaled \magstep3  % for chapter titles
\font\typf=cmcsc10                  % for authors
\font\typg=cmcsc10 scaled \magstep2 % for authors in titles
\global\def\sectitle#1\par{\bigbreak
  \leftline{\bf #1}
  \nobreak\medskip\vskip-\parskip
  \message{#1}
  \noindent}
\global\def\ssectitle#1\par{\bigbreak\medskip
  \leftline{\typc #1}
  \nobreak\bigskip\vskip-\parskip
  \message{#1}
  \noindent}
\def\theorem#1{\par\medbreak\noindent\bf Theorem\ #1.\enspace\sl\nobreak}
\def\endtheorem{\par\medskip\rm}

\def\proposition#1{\par\medbreak\noindent\bf
Proposition\ #1.\enspace\sl\nobreak}
\def\endproposition{\par\medskip\rm}
\def\corollary#1{\par\medbreak\noindent\bf
Corollary\ #1.\enspace\sl\nobreak}
\def\endcorollary{\par\medskip\rm}

\def\example#1{\par\medbreak\noindent\bf Example\ #1.\enspace\rm\nobreak}
\def\endexample{\par\medskip\rm}

\def\proof{\par\medbreak\noindent{\bf Proof.}\enspace}
\def\qed{\vbox{\hrule
  \hbox{\vrule\hbox to 5pt{\vbox to 8pt{\vfil}\hfil}\vrule}\hrule}}
\def\endproof{\unskip \nobreak \hskip0pt plus 1fill \qquad \qed% \par
\medskip\noindent}
\def\ref#1#2#3#4{{\item{[#1]}\typf #2 \it #3 \rm #4 \medskip}}
\def\R{\hbox{\rm I\kern-2pt R}}
\def\subR{_{\rm I\!R}}
\def\N{{\rm I\kern-2pt N}}
\def\Z{{\rm Z}\!\!{\bf Z}}
\def\C{\hbox{\rm \kern 3pt\raise .4pt
                          \hbox to 0.6pt{\vrule width .6pt height 6pt}
                          \kern-6.6pt C}}
\def\Q{\hbox{\rm \kern 3pt\raise .4pt
                          \hbox to 0.6pt{\vrule width .6pt height 6pt}
                          \kern-6.6pt Q}}
\def\subC{_{\rm 1\!\!\!C}}
\def\RP{\hbox{\rm I\kern-2pt R\kern-1.2pt P}}
\def\CP{\hbox{\C\kern-1.2pt\rm P}}
\def\Sum{\mathop{\textstyle{\sum}}\limits}
\def\lra{\longrightarrow}                % rightarrow
\def\iff{\Longleftrightarrow}

\def\sse{\subseteq}
\def\cir{{\circ}}

     % overline (high)

\def\es{\emptyset}
\def\sm{{\setminus}}
\def\ker{{\rm ker}}
\def\H{\widetilde{{\rm H}}}
\def\HH{{\rm H}}

\def\x{{\bf x}}

\def\eps{{\epsilon}}
\def\dd{{\rm d}}	%	differential form -d
\def\zb{\overline{z}}	%	z-bar
\def\BB{{\cal B}}
\def\codim{{\rm codim}}
\def\rank{{\rm rank}}
\def\sign{{\rm sign}}
\def\sing{{\rm sing}}
\def\BC{\hbox{\tt BC}}

\hrule
\bigskip
\centerline{\type On the Difference}
\medskip
\centerline{\type Between Real and Complex Arrangements}
\medskip

$$\vbox{ \halign{ # \quad & \hfil # \hfil & # \quad   \cr
& {\typg G\"unter M.~Ziegler}  & \cr
& Institut Mittag-Leffler      & \cr
& Aurav\"agen  17     	       & \cr
& S-18262 Djursholm, Sweden    & \cr }}$$

\medskip
\hrule

\medskip
\rightline{\sevenrm 3.~December 1991/e-print alg-geom/9202005}
\vskip2cm
\vfill

\item{\ }
{\bf Abstract. \ }
If $\BB$ is an arrangement of linear complex hyperplanes in $\C^d$, then
the following can be constructed from knowledge of its
intersection lattice:
\itemitem{(a)}
the cohomology groups of the complement [Br],
\itemitem{(b)}
the cohomology algebra of the complement [OS],
\itemitem{(c)}
the fundamental group of the complement, if $d\le 2$,
\itemitem{(d)}
the singularity link up to homeomorphism, if $d\le 3$,
\itemitem{(e)}
the singularity link up to homotopy type [ZZ].
\item{}
If $\BB'$ is, more generally, a $2$-arrangement in $\R^{2d}$
(an arrangement of real subspaces of codimension $2$ with even-dimensional
intersections), then the intersection lattice still determines
(a) the cohomology groups of the complement [GM] and (e) the homotopy type
of the singularity link [ZZ].
\item{}
We show, however, that for $2$-arrangements
the data (b), (c) and (d) are not determined by
the intersection lattice. They require the knowledge of extra
information on sign patterns, which can be computed as
determinants of linear relations, or (equivalently) as
linking coefficients in the sense of knot theory.
\eject

\ssectitle{1.\ Introduction}

Let $\BB=\{H_1,\ldots,H_n\}$ be an arrangement of
complex hyperplanes in $\C^d=\R^{2d}$.
We will only consider arrangements that are linear (the hyperplanes are
vector subspaces) and essential (the intersection of all the
hyperplanes is $\{0\}$).

The principal combinatorial structure associated with a complex
arrangement is the {\it intersection lattice}
$L_{\BB} :=\{\bigcap_{a\in A}H_a:A\sse \{1,\ldots,n\}\}$
of all intersections of
hyperplanes, ordered by reversed inclusion. This is a
geometric lattice (matroid), whose rank function is given by
complex codimension, $r(A)=\codim\subC(\bigcap_{a\in A}H_a)$.

Let $D_{\BB}:=S^{2d-1}\cap\bigcup\BB$ denote the {\it singularity link} and
let $C_{\BB}:=\C^d\sm\bigcup\BB$ denote the {\it complement} of the
arrangement. A by now classical result of Arnol'd, Brieskorn and
Orlik \& Solomon asserts that
a presentation of the cohomology algebra of $C_{\BB}$ can be constructed
from the data that are
encoded by the intersection lattice $L_{\BB}$, as follows.

\theorem{1.1} \ {\rm [A][Br][OS]} \
Let $\BB=\{H_1,\ldots,H_n\}$ be a complex arrangement in $\C^d$.
For every $H_a\in\BB$ choose a linear form $l_a\in(\C^d)^*$ that
defines it, such that $\ker(l_a)=H_a$ $(1\le a\le n)$. Then the
integral cohomology algebra of the complement is generated by the
classes
$$\omega_a\ :=\ {1\over 2\pi i}{\dd l_a\over l_a},$$
for $1\le a\le n$. It has a presentation of the form
$$0\ \lra\ I\ \lra\ \Lambda^*\Z^n
\ \ {\buildrel\pi\over\lra}\ \ \HH^*(C_{\BB} ;\Z)\ \lra\ 0,$$
defined by $\pi( e_a):=[\omega_a]$, where
$\{ e_1,\ldots,e_n\}$ denotes a basis of $\Z^n$.
The relation ideal $I$ is generated by the elements
$$\Sum_{i=0}^k (-1)^i\ e_{a_0}\wedge\ldots\wedge\widehat{\mathstrut e_{a_i}}
\wedge\ldots\wedge e_{a_k},$$
for circuits $A=\{a_0,\ldots,a_k\}$ of $L$, that is, for
the minimal subsets $A\sse\{1,\ldots,n\}$ with $r(A)<|A|$.
\endtheorem

Goresky \& MacPherson [GM, p.~257],
whose section title we have used for this paper, suggest to study the
following (seemingly) mild generalization.
A {\it $2$-arrangement} is a finite set $\BB'=\{H_1,\ldots,H_n\}$
of real vector subspaces of codimension $2$ in $\R^{2d}$ so that
every intersection $\bigcap_{a\in A}H_a$ has even codimension in $\R^{2d}$.
Again we assume $\bigcap\BB'=\{0\}$.
The combinatorial essence
of a ``complex structure'' can be studied by comparing the structure
of $2$-arrangements with that of complex arrangements.

The intersection lattice of a
$2$-arrangement is again a geometric lattice, where
real codimension corresponds to twice the lattice rank:
$2{\cdot}r(A)=\codim\subR(\bigcap_{a\in A}H_a)$.

The cohomology groups of the complements of $2$-arrangements were computed by
Goresky \& MacPherson using Stratified Morse Theory.
An alternative approach to the computation, via spectral sequences, is
provided by Vassiliev [Va] and Jewell, Orlik \& Shapiro [JOS].
A third proof, with homotopy methods, is given by
Ziegler \& \v Zivaljevi\'c [ZZ].
However, the algebra structure is not supplied by either approach.
The combinatorial method of Bj\"orner \& Ziegler [BZ] yields the following
information about it.

\theorem{1.1$'$} \ {\rm [GM][BZ]} \
Let $\BB=\{H_1,\ldots,H_n\}$ be a $2$-arrangement in $\R^{2d}$.
For every $H_a\in\BB'$ choose two linear forms $l_a,l_a'\in(\R^{2d})^*$ that
define it, such that $\ker(l_a)\cap\ker(l_a')=H_a$ $(1\le a\le n)$. Then the
integral cohomology algebra of the complement $C_{\BB}$ is generated by the
$1$-dimensional classes
$$\omega(l_a,l_a')\ :=\
{1\over 2\pi}{-l_a'\dd l_a+l_a\dd l_a'\over {l_a}^2+{l_a'}^2},$$
for $1\le a\le n$. It has a presentation of the form
$$0\ \lra\ I\ \lra\ \Lambda^*\Z^n
\ \ {\buildrel\pi\over\lra}\ \ \HH^*(C_{\BB};\Z)\ \lra\ 0,$$
defined by $\pi( e_a):=[\omega(l_a,l_a')]$, where
$\{e_1,\ldots,e_n\}$ denotes a basis of $\Z^n$.
The relation ideal $I$ is generated by elements of the form
$$\Sum_{i=0}^k \eps_i\cdot e_{a_0}\wedge\ldots\wedge
\widehat{\mathstrut e_{a_i}} \wedge\ldots\wedge e_{a_k},$$
for the circuits $A=\{a_0,\ldots,a_k\}$ of $L$,
with $\eps_i\in\{+1,-1\}$.
\endtheorem

In the following, we will show that the inability of Theorem $1.1'$
to determine the precise form of the presentation of $\HH^*(C_{\BB'};\Z)$
from the combinatorial data is not a weakness of
stratified Morse theory of [GM] and of the combinatorial set-up of [BZ].
In fact, the cohomology algebra $\HH^*(C_{\BB'};\Z)$,
and hence the homotopy type of the complement
of a $2$-arrangement, is not determined by the
combinatorial data!

\theorem{1.2}
There are two different $2$-arrangements $\BB$ and $\BB'$ of $2$-dimensional
linear subspaces in $\R^4$ whose intersection lattices coincide
(the corresponding matroid is the uniform matroid $U_{2,4}$),
but whose complements have non-isomorphic cohomology algebras.
\endtheorem

In the following section we will give an extensive analysis of the
topology of the $2$-arrangements of $4$
transversal $2$-subspaces in $\R^4$ (corresponding to
the uniform matroid $U_{2,4}$), and show how they
can sometimes be distinguished by the cohomology algebras
of their complements. In Section~3 the implications for the singularity links
of $2$-arrangements are derived. In Section~4 we obtain a
general method to compute a presentation of  $\HH^*(C_{\BB'};\Z)$
once equations for $\BB'$ are chosen.
Section~5 discusses the relation to the study of knots and links in $S^3$.

In the following we will denote $2$-arrangements by $\BB'$, and only
drop the prime in the case of a complex arrangement.

\ssectitle{2.\ Example}

In this section we consider $2$-arrangements in $\R^4$:
arrangements of $2$-dimensional
linear subspaces $\BB'=\{H_a:1\le a\le n\}$ in $\R^4$
that are {\it transversal}, that is, have pairwise intersection $\{0\}$.
They represent the uniform matroid $M=U_{2,n}$.
We will use coordinates $u,v,x,y$ on $\R^4$, which we
abbreviate as $w=u+iv$, $z=x+iy$ when the usual complex structure
(identification of $\C^2$ and $\R^4$) is chosen.

In suitable coordinates we can assume that
$$\eqalign{
H_1=&\ \{(w,z)\in\R^4: w=0\ \},\cr
H_2=&\ \{(w,z)\in\R^4: z=0\ \},\cr
H_3=&\ \{(w,z)\in\R^4: z=w\,\}.\cr}$$

We note here that the projection $\pi:\R^4\lra\R^2$, which maps
$(u,v,x,y)\longmapsto(u,v)$, $(w,z)\longmapsto w$,
makes the complement $C_{\BB'}$ into a fiber bundle over $\C^*$, whose
fiber is $\C$ minus $n{-}1$ points.
As a consequence of this we deduce that every $2$-arrangement in
$\R^4$ is a $K(\pi,1)$-arrangement: the long exact homotopy
sequence of the fiber bundle shows that the higher homotopy groups of
$C_{\BB'}$ vanish in this case.

The fiber bundle is trivial if $\BB$ is a complex arrangement [Or, Prop.~5.3]:
assume that the hyperplanes are $H_1=\{(w,z): w=0\}$
and $H_a=\{(z,w):z=\lambda_aw\}$ for $2\le a\le n$, then
$$\eqalign{\mu:\ C_{\BB}\ \ &\lra\ \
\C^*\times\C\sm\{\lambda_2,\ldots,\lambda_n\}\cr
(w,z)\ \ &\longmapsto\ \ (w,{z\over w})      \cr}$$
trivializes the bundle.
The bundle is usually non-trivial for $2$-arrangements in $\R^4$.
Our results of this section will imply that it is not in general
homotopy equivalent to a product space with a factor $\C^*$.

In the complex case, even more can be said.
For this, note that the singularity link
of a complex arrangement in $\C^2$ is a disjoint union of circles, and
each circle has a natural orientation, given by multiplication with
$e^{it}$.

\proposition{2.1}
Let $\BB_1$ and $\BB_2$ be two arrangements of $n$ hyperplanes
($1$-dimensional complex subspaces) in $\C^2$,
with singularity links $D_1$ and $D_2$.
Then every orientation-preserving homeo\-morphism $D_1\lra D_2$
can be extended to a homeomorphism $(S^3,D_1)\lra(S^3,D_2)$.
\endproposition

\proof
Assume that $\BB_i$ is given by $l_{i1}(w,z)=z$ and
$l_{ij}(w,z)=w-\lambda_{ij}z$
for $2\le j\le n$. Then any homeomorphism of the Riemann sphere
that fixes infinity and maps $\lambda_{1j}$ to $\lambda_{2j}$ for all
$j$ yields a homeo\-mor\-phism $(S^3,D_1)\lra(S^3,D_2)$. The fact that
the initial homeomorphism $D_1\lra D_2$ can be prescribed
arbitrarily now follows from surgery
along a tubular neighborhood of $D_1$ resp.\ $D_2$.
\endproof

Now we will restrict our attention to the case $n=4$.
It is not hard to see [VD, p.~1038]
that there are three isotopy classes of arrangements.
One class contains the complex arrangements, which we denote by $\BB$.
The image of a complex arrangement $\BB$ after a reflection in $\R^4$
is not isotopic to a complex arrangement, but clearly isomorphic as
a $2$-arrangement. The third class is represented by the
$2$-arrangement $\BB'$ in the second case of the following example.

\example{2.2}
As the {\bf first case} we consider the arrangement $\BB$
$$\BB:\ \ \left\{\matrix{
H_1=&\{(w,z)\in\R^4: &w=0\}	\cr
H_2=&\{(w,z)\in\R^4: &z=0\}	\cr
H_3=&\{(w,z)\in\R^4: &z=w\}	\cr
H_4=&\{(w,z)\in\R^4: &z=2w\}. 	\cr}\right.
$$
This is a complex (in fact: complexified real) arrangement.
Its cohomology algebra has by Theorem 1.1 the following
presentation:
$$\HH^*(C_{\BB};\Z)\ \cong\ \Lambda^*\Z^4\big/
\left\langle\matrix{
+e_{12}-e_{13}+e_{23}\cr
+e_{12}-e_{14}+e_{24}\cr
+e_{13}-e_{14}+e_{34}\cr
+e_{23}-e_{24}+e_{34}\cr}\right\rangle $$
where the last relation is a consequence of the first three.
\medskip

As the {\bf second case} we consider the arrangement $\BB'$
$$\BB':\ \ \left\{\matrix{
H_1=&\{(w,z)\in\R^4: &w=0	\}\	\cr
H_2=&\{(w,z)\in\R^4: &z=0	\}\	\cr
H_3=&\{(w,z)\in\R^4: &z=w	\}\	\cr
H_4=&\{(w,z)\in\R^4: &z=2\overline w\}. \cr}\right.
$$
This arrangement is not linearly isomorphic to a complex one.
However, its cohomology algebra has by Theorem $1.1'$
a very similar presentation. With the method of Theorem 4.1 below, one can
determine the signs in the presentation:
$$\HH^*(C_{\BB'};\Z)\ \cong\ \Lambda^*\Z^4\big/
\left\langle\matrix{
+e_{12}-e_{13}+e_{23}\cr
+e_{12}+e_{14}+e_{24}\cr
-e_{13}-e_{14}+e_{34}\cr
+e_{23}-e_{24}-e_{34}\cr}\right\rangle $$
where the last relation is a consequence of the first three.
\endexample

In both cases the broken circuit complex
$$\BC(U_{2,4})=\{\es,1,2,3,4,12,13,14\}$$
indexes a basis of the cohomology module, that is, the (classes of)
$e_1,e_2,e_3,e_4$ induce a $\Z$-basis of $\HH^1$, while
$e_{12},e_{13},e_{14}$ induce a $\Z$-basis of $\HH^2$, and
$\HH^3=\HH^4=0$, see [BZ, Sect.~7].
In particular, the cohomology modules
$\HH^*(C_{\BB};\Z)$ and $\HH^*(C_{\BB'};\Z)$ are linearly isomorphic.
Their difference is hidden in the multiplicative structure.
For the following theorem we do not assume that
an isomorphism maps generators to generators, but argue with
an invariant construction. It was inspired by [F1], although
Falk's invariants do not suffice to distinguish the
algebras $\HH^*(C_{\BB};\Z)$ and $\HH^*(C_{\BB'};\Z)$.

\theorem{2.3}
The cohomology algebras $\HH^*(C_{\BB};\Z)$ and $\HH^*(C_{\BB'};\Z)$
are not isomorphic as graded $\Z$-algebras.
\endtheorem

\proof
Let $A$ denote any of the two algebras and let $A^1$ be its
$1$-dimensional part. Then $A$ has a presentation of the
form
$$0\ \lra\ I\ \lra\ \Lambda^*A^1\ \lra\ A\ \lra\ 0,$$
where $I$ is again a graded ideal. Here $I^1=0$ by
construction, while $I^2$ has rank $3$.
We consider the map
$$\kappa:\ I^2\ \otimes I^2\ \ \lra\ \ \Lambda^4A^1$$
induced by multiplication in $\Lambda^*A^1$.
In the case we are considering $A^1\cong\Z^4$, so
$\Lambda^4A^1\cong\Z$, and $\kappa$ defines a symmetric
bilinear form on $I^2$.

A direct calculation shows that $\kappa$ vanishes identically
for $\HH^*(C_{\BB};\Z)$. This can also be derived from the
K\"unneth formula, since $C_{\BB}$ is a product space.

However, for $\BB'$ the bilinear form $\kappa$ has rank $2$:
with respect to the basis $\{e_{12}-e_{13}+e_{23},$
$e_{12}+e_{14}+e_{24}, -e_{13}-e_{14}+e_{34}\}$ it is represented by
the matrix
$$\pmatrix{0 & 2& 0\cr
	   2 & 0&-2\cr
	   0 &-2& 0\cr}.$$
This proves $\HH^*(C_{\BB};\Q)\not\cong \HH^*(C_{\BB'};\Q)$.
\endproof

\ssectitle{3.\ Links of 2-arrangements}

The singularity link of any $2$-arrangement is
homotopy equivalent to a wedge of spheres if $d\ge 3$
[BZ, Thm.~6.6] [ZZ, Cor.~3.3] and it is a disjoint union of
circles if $d=2$.

In fact, there is a certain plausibility to
the conjecture that for complex arrangements the
singularity links are determined
up to homeomorphism by the intersection lattices,
for all $d$. This is a very strong conjecture, which would
(with the ideas below) imply the same for the complements
of complex arrangements, which is much stronger than the
notorious conjecture [Or] for the homotopy type.
In this section, we will prove this fact for $d\le 3$, and
then use the example of Section 2 to disprove a similar
statement in the case of $2$-arrangements.

\theorem{3.1}
The intersection lattice determines the
singularity link of a complex arrangement in $\C^3$ up to homeomorphism.
\endtheorem

\proof
Let $\BB$ be a complex arrangement of $n$ $2$-dimensional
subspaces in $\C^3$, with intersection lattice $L$. The singular
set of $D$ is a disjoint union of $k$ circles, where
$k$ is the number of coatoms (elements of rank $2$) in $L$.
Thus we construct $\sing(D)$ as a set of $k$ disjoint oriented circles, where
the orientation is supposed to be the natural one corresponding to
multiplication with $e^{it}$. Now we glue $n$ $3$-spheres into the
given set of oriented circles. The attaching maps exist and are unique by
Proposition 2.1.
\endproof
\eject

\theorem{3.2}
The intersection lattice does not determine the
singularity link of a $2$-arrangement in $\R^6$ up to homeomorphism.
\endtheorem

\proof
We consider generic $2$-arrangements
$\hat\BB'=\{H_1,H_2,H_3,H_4,H_5\}$ representing the matroid $U_{3,5}$,
that is, arrangements of five $4$-dimensional subspaces in $\R^6$ so that
the intersection of any three of them is $\{0\}$.
Their singularity links are unions of five copies $S_i:=H_i\cap S^5$
of $S^3$, pairwise intersecting in circles.
We notice that the non-singular parts
$S_i^\cir:=S_i\sm\bigcup_{j\neq i} S_j$ are easily identified by
local cohomology. Each of these parts $S_i^\cir$ is homeomorphic
to the complement of the restriction
$\hat\BB'|{H_i}$, which is a $2$-arrangement of four
$2$-subspaces in $\R^4$,
as discussed in Section~2. If $\hat\BB'$ is a complex
arrangement, then the restrictions $\hat\BB'|H_i$ are complex as well.

However, for $\hat\BB'$ given by the equations
$$\eqalign{
H_1:&\ \  z_1=0,			\cr
H_2:&\ \  z_2=0,			\cr
H_3:&\ \  z_3=0,			\cr
H_4:&\ \  z_1-z_2+z_3=0,		\cr
H_5:&\ \  z_1-2\zb_2+3z_3=0,		\cr}$$
we find that $\hat\BB'|H_3$ is isomorphic to the arrangement
$\BB'$ considered in
Example~2.2, so $S_3^\cir$ is homeomorphic to $C_{\BB'}$, and hence it is
not homeomorphic to a non-singular part of the singularity link of a complex
arrangement.
\endproof

\ssectitle{4.\ Cohomology of 2-arrangements}

In this section we describe a method to compute the relations in the
cohomology algebra of any $2$-arrangement.
It relies on the representation of cohomology
classes by the corres\-ponding differential forms of real deRham theory,
and it exploits the passage to complex deRham theory in the case of
subarrangements that have a complex structure, like those corresponding
to the circuits of the matroid. It seems desirable to derive
a presentation in the combinatorial framework and generality of [BZ];
however, this has not yet been achieved.

We will need the relation between the real and the complex
differential form representing the cohomology class of a
complex hyperplane. For this assume that coordinates have been chosen so that
the hyperplane $H$ is represented by $z=0$, which with $z=x{+}iy$
corresponds to real equations $x=y=0$. Then straightforward
computations show that
$${1\over 2\pi i}\left({\dd z\over z}+{\dd \zb\over \zb}\right)\ =\
  {1\over 2\pi i}\dd\log(x^2+y^2),$$
which is an exact form, while
$${1\over 2\pi i}\left({\dd z\over z}-{\dd \zb\over \zb}\right)\ =\
  {1\over \pi}{-y\dd x+x\dd y\over x^2+y^2},$$
which is twice the real differential form that
represents the cohomology class of $\C^d\sm H$. Thus
$$\left[{ 1\over 2\pi i}{\dd  z \over  z }\right]\ =\
  \left[{-1\over 2\pi i}{\dd \zb\over \zb}\right]\ =\
  \left[{ 1\over 2\pi}{-y\dd x+x\dd y\over x^2+y^2}\right].$$

\theorem{4.1}
Let a $2$-arrangement $\BB'=\{H_1,\ldots,H_n\}$ in $\R^{2d}$ be
given by
$$H_a=\{\x\in\R^{2d}: l_a(\x)=l_a'(\x)=0\},$$
where the $l_a,l_a':\R^{2d}\lra\R$ are linear forms so that
$$\rank\{l_a,l_a':a\in A\}\ =\ \cases{
\hbox{\rm even}	& for all $A\sse\{1,\ldots,n\}$,\cr
2		& for all $A=\{a\}$, 		\cr
4		& for all $A=\{a,b\}$, 		\cr
2d		& for $A=\{1,\ldots,n\}$.	\cr}$$
To every $H_a\in\BB'$ associate the differential form
$$\omega(l_a,l_a')\ :=\
{1\over 2\pi}{-l_a'\dd l_a+l_a\dd l_a'\over {l_a}^2+{l_a'}^2},$$
which is a closed form on $\R^{2d}\sm H_a$ that is normalized to have residue
$\pm1$.
The relations between the corresponding cohomology classes
can be constructed as follows. Let $A=\{a_0,a_1,\ldots,a_k\}$
be a circuit of $L_{\BB'}$, so
there are two real linear dependencies of the form
$$\eqalign{
\Sum_{j=0}^k\ \alpha_j l_{a_j} +\beta_j l'_{a_j}\ &=\ 0,\cr
\Sum_{j=0}^k\ \gamma_j l_{a_j} +\delta_jl'_{a_j}\ &=\ 0,\cr}$$
with $\alpha_0=\delta_0=-1$, $\beta_0=\gamma_0=0$.
These induce the relation
$$\Sum_{j=0}^k (-1)^j
\sign\left|\matrix{\alpha_j&\beta_j\cr\gamma_j&\delta_j\cr}\right|
\ \omega(l_{a_1},l_{a_1}')\wedge\ldots\wedge
\widehat{\mathstrut\omega(l_{a_j},l_{a_j}')}
\wedge\ldots\wedge\omega(l_{a_k},l_{a_k}')\ \ \sim\ \ 0$$
in the cohomology algebra $\HH^*(C_{\BB'};\Z)$.
\endtheorem

\proof
The conditions on the forms $l_a,l_a'$ assure that they define
a $2$-arrangement.
The differential forms $\omega(l_a,l_a')$
generate $\HH^*(C_{\BB'};\Z)$, by Theorem $1.1'$.

To derive the relations we
construct coordinates $x_j,y_j$ for $\R^{2d}$ so that
$$\eqalign{
x_j\ =&\ \alpha_j l_{a_j} +\beta_j l'_{a_j}\cr
y_j\ =&\ \gamma_j l_{a_j} +\delta_jl'_{a_j}\cr}$$
for $1\le j\le k$ ---
this is possible since $A$ is a circuit, so $\{l_a,l'_a:a\in A\sm a_0\}$
is linearly independent. The even rank condition furthermore guarantees
$\left|\matrix{\alpha_j&\beta_j\cr\gamma_j&\delta_j\cr}\right|\neq0$.
Now we observe, by computing the residues, that
$$\omega(x_{a_j},y_{a_j})\ \sim\
\sign\left|\matrix{\alpha_j&\beta_j\cr\gamma_j&\delta_j\cr}\right|
\cdot\omega(l_{a_j},l_{a_j}').$$
Introducing complex coordinates $z_j:=x_j{+}iy_j$, we get
that $H_{a_0}$ has the (complex!) equation $z_1{+}\ldots{+}z_k=0$,
while $H_{a_j}$ is given by $z_j=0$ for $1\le j\le k$,
and thus
$$\omega_j
\ :=\ {1\over 2\pi i}{\dd z_j\over z_j}
\ \sim\ \omega(x_j,y_j)
\ \sim\ \sign\left|\matrix{\alpha_j&\beta_j\cr\gamma_j&\delta_j\cr}\right|\cdot
\omega(l_{a_j},l_{a_j}').$$
Thus the relation
$$\Sum_{j=0}^k (-1)^j\ \omega_0\wedge\ldots\wedge\widehat{\mathstrut\omega_i}
\wedge\ldots\wedge\omega_k\ =\ 0$$
for the complex arrangement $\{H_{a_0},H_{a_1},\ldots,H_{a_k}\}$
translates into the desired formula.
\endproof

Note that the formula of Theorem 4.1 specializes to the
Orlik-Solomon relations in the case of a complex arrangement:
for a complex arrangement we can write the defining forms
as $l_a+il_a'$, and the relation corresponding to a circuit takes the
form
$$\Sum_{j=0}^k(\alpha_j+i\beta_j)(l_{a_j}+il_{a_j}')=0\qquad\iff\qquad
\left\{\matrix{
\Sum_{j=0}^k\hphantom{+}\alpha_jl_{a_j}+\beta_j l_{a_j}'&=0\cr
\Sum_{j=0}^k          - \beta_j l_{a_j}+\alpha_jl_{a_j}'&=0\cr}\right.$$
Thus for the formula of Theorem 4.1 we get the special case
$\gamma_j=-\beta_j$ and $\delta_j=\alpha_j$, so that \
$\sign\left|\matrix{\alpha_j&\beta_j\cr\gamma_j&\delta_j\cr}\right|=
\sign(\alpha_j^2+\beta_j^2)=+1$.

\ssectitle{5.\ Link invariants}

The classification of $2$-arrangements in $\R^4$ is clearly equivalent to the
study of
\item{--}
arrangements of disjoint great circles in $S^3$,
\item{--}
arrangements of skew lines in $\RP^3$,
\item{--}
arrangements of affine skew lines in $\R^3$,

\noindent
as is e.g.\ stressed by Viro [V] and by
Viro \& Drobotukhina [VD, p.~1046]. The corresponding equivalence relation
on line arrangements is there called {\it rigid isotopy}.

Considering arrangements of circles in $S^3$ as {\it links}, one is lead to
study to what extent link invariants can distinguish
equivalence classes of $2$-arrangements in $\R^4$.
In particular, the $(2\times2)$-determinants derived in Section~4, which
determine the sign pattern of the relations in cohomology, are just
linking numbers of the corresponding (oriented) circles. In the
description as links in the $3$-sphere it is well-known [VD, p.~1034] that
for every triple of circles we get a {\it linking coefficient} $\pm1$
that does not depend on the order or the chosen orientations.
These linking coefficients of triples
are sufficient to distinguish the arrangements $\BB$ and $\BB'$ of Example~2.2.
We refer to [V], [VD] and [M] for this approach.

The key problem here is that links do not in general
determine the homotopy type of their complements in $S^3$
(see e.g.\ [R, p.~62]), although this might be true for the
special type of links that correspond to arrangements.
Therefore the results of [V], [VD] and [M] do not immediately distinguish
complements of $2$-arrangements.
In fact, the relation between the ``new'' link invariants and
``classical'' data like the fundamental group of the
complement is still obscure [Bi, p.~59].

A presentation of the fundamental group of the complement of a $2$-arrangement
can easily be computed --- the standard method due to
Wirtinger (see [R]) derives it from a planar projection; since the
links we consider are closed braids, an equivalent (but more systematic) way
is given by Artin's approach [Ar]. However, even for the simple case of
the $2$-arrangement of four subspaces in $\R^4$ the corresponding
links have projections with $12$ crossings, so these methods become
unwieldy. From the description as a fiber bundle, one sees that
the fundamental group $\pi$ in the case of a complex arrangement is a
product of $\Z$ with the free group $F\langle t_1,t_2,t_3\rangle$
on three generators.
In the case of the arrangement $\BB'$, we find that $\pi'$ is a non-trivial
solution of the extension problem
$$F\langle t_1,t_2,t_3\rangle\ \lra\ \pi'\ \lra\ \Z\lra0.$$
However, there seems to be no simple or direct way to describe the
homotopy group $\pi'$.

\corollary{5.1}
The fundamental groups $\pi,\pi'$ of the complements of the arrangements
$\BB$ and $\BB'$ of Example 2.2 are not isomorphic.
\endcorollary

\proof
We have seen that $C_{\BB}\simeq K(\pi,1)$ and $C_{\BB'}\simeq K(\pi',1)$
in Section~2, and that these spaces have non-isomorphic cohomology algebras
in Theorem~2.3. Hence $\pi\not\cong\pi'$.
\endproof

\example{5.2} \ {\rm [M, Ass.~3] [VD, p.~1043]} \
There are two $2$-arrangements $\BB'$ and $\BB''$ of six two-dimensional
transversal subspaces in $\R^4$ with the following properties:
\item{--}
the cohomology algebras $\H^*(C_{\BB'};\Z)$ and $\H^*(C_{\BB''};\Z)$
are isomorphic, because the subspaces in the arrangements can be labeled
and oriented in such a way that the pairwise linking numbers coincide,
\item{--}
the pairs $(S^3,D')$ and $(S^3,D'')$ are not homeomorphic,
since they represent inequivalent links in $S^3$ that can be distinguished by
link polynomials.

\noindent
We do not know whether the complements $S^3\sm D''\simeq C_{\BB''}$ and
$S^3\sm D'\simeq C_{\BB'}$ are homotopy equivalent
or, equivalently (by the argument of Corollary 5.1), whether
the fundamental groups coincide.
\endexample

More generally, we do not know whether
the complements of two $2$-arrangements must be homotopy equivalent once
their cohomology algebras are isomorphic.

It seems that the analysis of the cohomology algebra (as in Theorem 2.3)
is simpler than any computation of the fundamental group of a
$2$-arrangement. However, we are not aware of any systematic study of the
cohomology algebra of the complement of a link (compare e.g.\ [R, p.~50]).
The linear structure of this algebra is
determined by the number of components of the link,
because of Alexander duality. But, as we have seen
in Section 2, the multiplicative structure encodes
non-trivial information.

However, we note that the complement of every $2$-arrangement is
{\it formal} in the sense of rational homotopy theory [GrH, p.~158].
In fact, by Theorem 1.1$'$ the cohomology algebra $\H^*(C_{\BB'};\Z)$
can be represented by a subalgebra of the real deRham complex on $C_{\BB'}$.
With the argument of [F2, p.~546] this implies that $C_{\BB'}$ is a formal
space. In particular, there is no ``higher order'' cohomology information
(like the Massey products used in [GrM, Sect.~IIIX.C]) contained in the
real deRham complex of $C_{\BB'}$.
\vfill

\sectitle{Acknowledgements.}

I am grateful to Anders Bj\"orner, Michael Falk, Rade \v Zivaljevi\'c and
in particular to Boris Shapiro for many useful discussions,
explanations and references.
\eject

\ssectitle{References.}

\ref{A}
{V.I.~Arnol'd:}
{The cohomology ring of the colored braid group,}
{{\sl Mathematical Notes} {\bf 5} (1969), 138-140.}

\ref{Ar}
{E.~Artin:}
{Theorie der Z\"opfe,}
{{\sl Abhandlungen aus dem Mathematischen
Seminar der Hamburgischen Universit\"at} {\bf 4} (1926), 47-72.}

\ref{Bi}
{J.~S.~Birman:}
{Recent developments in braid and link theory,}
{{\sl Math.\ Intelligencer} {\bf 13} (1991), 52-60.}

\ref{BZ}
{A.~Bj\"orner \& G.~M.~Ziegler:}
{Combinatorial stratification of complex arrangements,}
{{\sl Journal Amer.\ Math.\ Soc.} {\bf 5} (1992), in press.}

\ref{Br}
{E.~Brieskorn:}
{Sur le groupe de tresses (d'apr\` es V.~I.~Arnol'd),}
{S\'eminaire Bourbaki 24e ann\'ee 1971/72, Lecture Notes in
Mathematics {\bf 317} (1973), 21-44.}

\ref{F1}
{M.~J.~Falk:}
{On the algebra associated with a geometric lattice,}
{{\sl Advances in Math.} {\bf 80} (1990), 152-163.}

\ref{F2}
{M.~J.~Falk:}
{The minimal model of the complement of an arrangement of hyperplanes,}
{{\sl Transactions Amer.\ Math.\ Soc.} {\bf 309}, 543-556.}

\ref{GM}
{M.~Goresky \& R.~MacPherson:}
{Stratified Morse Theory,}
{{\sl Ergebnisse der Mathematik und ihrer Grenzgebiete},
3.~Folge, Band~14, Springer 1988.}

\ref{GrM}
{P.~A.~Griffiths \& J.~W.~Morgan:}
{Rational Homotopy Theory and Differential Forms,}
{Pro\-gress in Mathematics {\bf 16}, Birkh\"auser 1981.}

\ref{JOS}
{K.~Jewell, P.~Orlik \& B.~Z.~Shapiro:}
{On the cohomology of complements to arrangements of affine subspaces,}
{preprint 1991.}

\ref{M}
{V.~F.~Mazurovskii:}
{Kauffmann polynomials of non-singular configurations of projective lines,}
{{\sl Russian Math.\ Surveys} {\bf 44} (1989), 212-213.}

\ref{Or}
{P.~Orlik:}
{Introduction to Arrangements,}
{{\sl CBMS Regional Conference Series in Mathematics} {\bf 72},
Amer.\ Math.\ Soc., Providence RI, 1989.}

\ref{OS}
{P.~Orlik \& L.~Solomon:}
{Combinatorics and topology of complements of hyperplanes,}
{{\sl Inventiones math.} {\bf 56} (1980), 167-189.}

\ref{R}
{D.~Rolfsen:}
{Knots and Links,}
{Mathematics Lecture Notes Series {\bf 7}, Publish or Perish,
Berkeley CA, 1976.}

\ref{Va}
{V.~A.~Vassiliev:}
{Complements to discriminants of smooth mappings,}
{conference notes (International Centre for Theoretical Physics, Trieste),
1991.}

\ref{V}
{O.~Ya.~Viro:}
{Topological problems concerning lines and points of three-dimensional space,}
{{\sl Soviet Math.\ Dokl.} {\bf 32} (1985), 528-531.}

\ref{VD}
{O.~Ya.~Viro \& Yu.~V.~Drobotukhina:}
{Configurations of skew lines,}
{{\sl Leningrad J.~Math.} {\bf 1} (1990), 1027-1050.}

\ref{Z}
{G.~M.~Ziegler:}
{Combinatorial Models for Subspace Arrangements,}
{Habilitations\-schrift, TU~Berlin 1992, in preparation.}

\ref{ZZ}
{G.~M.~Ziegler \& R.~T.~\v Zivaljevi\'c:}
{Homotopy types of arrangements via diagrams of spaces,}
{Report No. \ (1991/92), Institut Mittag-Leffler, December 1991.}
\vfill

\noindent
{\bf Address after April 1, 1992:}
Konrad-Zuse-Zentrum f\"ur
Informationstechnik Berlin,
Heilbronner Str.~10,
W-1000 Berlin 31,
Germany
\medskip

\noindent
{\bf E-mail:}
{\tt ziegler@ml.kva.se} \ or \ {\tt ziegler@zib-berlin.de}
\eject
\end